\documentclass[12pt]{amsart}
\DeclareFontFamily{OML}{script}{}
\DeclareFontShape{OML}{script}{m}{it}
{ <5-20> rsfs10 }{}
\DeclareMathAlphabet{\mathscript}{OML}{script}{m}{it}

\renewcommand{\mathcal}[1]{{\mathscript #1}\hspace{0.2ex}}
\usepackage{calligra}
\usepackage{pstricks}
\usepackage{color}
\usepackage{float}
\ifx\red\undefined
\newcommand{\red}{\color{red}}

\fi

\usepackage{graphicx,graphics}
\usepackage{cite}
\usepackage{pifont}
\usepackage{amsthm}

\usepackage{amscd}
\usepackage{amsmath}
\usepackage{latexsym}
\usepackage{amsfonts}
\usepackage{amssymb}
\usepackage{color}
\usepackage[ansinew]{inputenc}
\usepackage[encapsulated]{CJK}

\ifx\text\undefined
\newcommand{\text}{\mbox}
\fi
\ifx\operatorname\undefined
\newcommand{\operatorname}{\mathop}
\fi

\newcommand\be{\begin{equation}}
\newcommand\ee{\end{equation}}
\newcommand\bea{\begin{eqnarray}}
\newcommand\eea{\end{eqnarray}}
\newcommand\beaa{\begin{eqnarray*}}
\newcommand\eeaa{\end{eqnarray*}}

\newenvironment{eqa}{\begin{equation}%
  \begin{array}{rcl}}{\end{array}\end{equation}}
\newcommand\beqa{\begin{eqa}}
\newcommand\eeqa{\end{eqa}}

\numberwithin{equation}{section}

\renewcommand{\tilde}{\widetilde}

\newtheorem{thm}{Theorem}[section]

\newtheorem{lem}{Lemma}[section]

\newtheorem{rem}{Remark}[section]

\newcommand{\void}[1]{}

\numberwithin{equation}{section}

\begin{document}\begin{CJK}{UTF8}{gkai}

\title[Asymptotic stability]{Asymptotic stability for a free boundary tumor model with a periodic supply of external nutrients}
\author[Huang]{Yaodan Huang}

\address[1]{School of Mathematics, Sun Yat-Sen University, Guangzhou, 510275, China}
\email{huangyd35@mail.sysu.edu.cn}
\thanks{Corresponding author: Yaodan Huang}

\maketitle

\begin{abstract}
For tumor growth, the morphological instability provides a mechanism for invasion via tumor fingering and fragmentation. This work considers the asymptotic stability of a free boundary tumor model with a periodic supply of external nutrients. The model consists of two elliptic equations describing the concentration of nutrients and the distribution of the internal pressure in the tumor tissue, respectively. The effect of the parameter $\mu$ representing a measure of mitosis on the morphological stability are taken into account. It was recently established in \cite{HZH3} that there exists a critical value $\mu_\ast$ such that the unique spherical periodic positive solution is linearly stable for $0<\mu<\mu_\ast$ and linearly unstable for $\mu>\mu_\ast$. In this paper, we further prove that the spherical periodic positive solution is  asymptotically 
stable for $0<\mu<\mu_\ast$ for the fully nonlinear problem. 
\end{abstract}

\section{Introduction}
Within the past several decades a number of mathematical models have been developed and studied that aimed at describing the evolution to carcinomas in the form of free boundary problems. Choices have to be made for a variety of different models to explore the mechanisms of tumor growth, such as several types of cells (see, for instance, \cite{WK2,TP}), the impact of inhibitors (cf. \cite{BC1}). Accordingly, many significant results have been obtained for rigorous mathematical analysis and numerical simulation of such tumor models; see \cite{CLN,C1,CE2,EM3,FL,HHH2,HHH3,HZH1,HZH2,HZH4,LL,SHZ,W0,WC2,WZ1,WZ,ZH,ZH2,ZX,ZC} and the references given there. Lowengrub et al. \cite{LFJC} provided a systematic survey of tumor model studies.


Recent study indicates that human beings and many animals have regular living and feeding activities, related to the biological rhythm,
the concentration of nutrients in their blood may change periodically over time \cite{FD}. Motivated by this study, in this paper, we shall consider a free boundary tumor model describing the growth of tumors with a periodic supply of external nutrients:
\begin{eqnarray}
   \Delta\sigma & =& \sigma \ \ \ \ \ \qquad \  \ \  \ \mbox{in} \ \Omega(t),\ \ t>0,\label{1.2} \\
  -\Delta p &=& \mu(\sigma-\tilde{\sigma}) \ \ \ \ \ \mbox{in}\ \Omega(t),\ \ t>0, \label{1.4}\\
  \sigma &=& \phi(t) \ \ \ \qquad \ \ \mbox{on}\ \partial\Omega(t),\ \ t>0,\label{1.3} \\
  p &=& \gamma\kappa \ \ \ \ \qquad \ \ \ \mbox{on}\ \partial\Omega(t),\ \ t>0,\label{1.5}\\
  V_n&=&-\frac{\partial p}{\partial \vec{n}}  \ \qquad \ \ \ \mbox{on}\ \partial\Omega(t),\ \ t>0.\label{1.6}
\end{eqnarray}
Here $\Omega(t)\subseteq\mathbb{R}^3$ is an apriorily unknown bounded domain occupied by the tumor at time $t$, $\sigma=\sigma(x,t)$ and $p=p(x,t)$ are unknown functions representing the nutrient concentration in the tumor region and the pressure between tumor cells, respectively. $V_n$ is the velocity of the free boundary in the direction $\vec{n}$ which is the unit outward normal. $\mu$, $\widetilde{\sigma}$ and $\gamma$ are positive constants, among which $\mu$ measures the aggressiveness of the tumor, $\widetilde{\sigma}$ is the threshold nutrient concentration for apoptosis, and $\gamma$ represents the surface tension coefficient of the tumor surface $\partial\Omega(t)$. $\kappa$ is the mean curvature of $\partial\Omega(t)$ whose sign is defined by the convention that convex hypersurfaces are associated with 
positive mean curvature; $\kappa=1/R(t)$ if $\Omega(t)$ is a ball of radius $R(t)$. Besides, $\phi(t)$ accounts for the external concentration of nutrients which we shall assume is a smooth, positive, periodic function with period $T$.

The problem (\ref{1.2})-(\ref{1.6}) is a natural extension of that proposed by Bryne and Chaplain in \cite{BC1}: the action of a periodic supply of external nutrients is considered, so that (\ref{1.3}) holds. In this model, the quasi-stationary approximation of the 
diffusion equation (\ref{1.2}) combined 
with the boundary condition (\ref{1.3}) are used to describe the distribution of the nutrient concentration. Accordingly, the problem (\ref{1.2})-(\ref{1.6}) is called the \emph{quasi-stationary model}. In addition, assuming that the extracellular matrix is a porous medium where Darcy's law $\vec{V}=-\nabla p$ holds, the law of conservation of mass $\mbox{div}\vec{V}=\mu(\sigma-\widetilde{\sigma})$ yields the equation (\ref{1.4}). Equation (\ref{1.6}) is the classical Stefan condition for free boundary $\partial\Omega(t)$. Finally, the cell-to-cell adhesiveness leads to the relation (\ref{1.5}) (see \cite{B,BC3,G2}). If the stationary diffusion equation (\ref{1.2}) is replaced by its non-stationary version
\begin{equation*}\label{0.67}
  c\sigma_t=\Delta\sigma-\sigma  \qquad \qquad \ \mbox{in} \ \Omega(t),\ \ t>0,
\end{equation*}
where the positive constant $c=T_{\mbox{diffusion}}/T_{\mbox{growth}}$ is the ratio of the nutrient diffusion time scale, $T_{\mbox{diffusion}}\approx1\mbox{min}$, to the tumor-cell doubling time scale, $T_{\mbox{growth}}\approx1\mbox{day}$, so that $c\ll1$ (see \cite{BC1}), the corresponding problem of (\ref{1.2})-(\ref{1.6}) is called the \emph{fully non-stationary model}.

If the external concentration of nutrients is assumed to be a constant, i.e., the boundary condition (\ref{1.3}) is reduced to $\sigma=\overline{\sigma}=\mbox{const.}$, rigorous mathematical analysis including existence and stability theorems have been established. For the fully non-stationary model, Friedman and Reitich \cite{FR1} proved that under the assumption $0<\widetilde{\sigma}<\overline{\sigma}$, there exists a unique radially symmetric stationary solution $(\sigma_S(r),p_S(r),R_S)$. They also proved that there exists a family of symmetric-breaking bifurcation branches of stationary solutions bifurcating from this unique radial stationary solution $(\sigma_S(r),p_S(r),R_S)$ in \cite{FR3}.
In the sequel, it was proved in \cite{BF} that $(\sigma_S(r),p_S(r),R_S)$ is asymptotically stable with respect to non-radial perturbations for $\mu$ sufficiently small. This work was refined by Friedman and Hu \cite{FH1,FH2}, that is, they determined the threshold value $\mu^\ast$ such that for $0<\mu<\mu^\ast$ the trivial solution is asymptotically stable under small non-radial perturbations, while for $\mu>\mu^\ast$ the trivial solution is unstable. Later on, Cui and Escher extended the asymptotic stability result to the general case that the nutrient consumption rate and the tumor cell proliferate rate are general increasing functions. They found a positive critical value $\gamma^\ast$ (surface tension coefficient in (\ref{1.5})) for which the radial stationary solution changes from instability to stability under non-radial perturbations for the quasi-stationary model in \cite{CE1} and for the fully non-stationary model with small $c$ in \cite{C3}, respectively.

In this paper, we study the quasi-stationary model (\ref{1.2})-(\ref{1.6}) with the nutrient supply given by a periodic function $\phi(t)$. We expect to determine how the tumor will evolve as time goes  to infinity under the effect of this periodic external nutrient supply. To give a precise statement of our main result, let us make some preparations. By a rescaling if necessary, we take $\gamma=1$. It is not difficult to verify (see also \cite{HX}) that the radially symmetric $T$-periodic positive solution of the problem (\ref{1.2})-(\ref{1.6}), which we denote by $(\sigma_\ast(r,t),p_\ast(r,t),R_\ast(t))$,
 is of the form
\begin{equation}\label{2.14}
\sigma_\ast(r,t)=\phi(t)\frac{R_\ast^{1/2}(t)}{I_{1/2}(R_\ast(t))}\frac{I_{1/2}(r)}{r^{1/2}}, \ \ \ 0<r<R_\ast(t),
\end{equation}
\begin{equation}\label{2.15}
p_\ast(r,t)=-\mu\sigma_\ast(r,t)+\frac{1}{6}\mu\tilde{\sigma}r^{2}+\frac{1}{R_\ast(t)}+\mu\phi(t)-\frac{1}{6}\mu\tilde{\sigma}R_\ast^{2}(t), \ \ \ 0<r<R_\ast(t),
\end{equation}
and $R_\ast(t)$ satisfies
\begin{equation}\label{2.16}
\begin{split}
  \frac{dR_\ast(t)}{dt}=\mu R_\ast(t)\Big\{\phi(t)P_0(R_\ast(t))-\frac{\widetilde{\sigma}}{3}\Big\},
\end{split}
\end{equation}
where $I_n(r)$ is the modified Bessel function of order $n$ and $P_0(r)$ is defined by (\ref{2.1}).
As was analyzed in \cite{WX}, for $\frac{1}{T}\int_0^T\phi(t)dt>\widetilde{\sigma}$, there exists a unique $T$-periodic positive solution $R_\ast(t)$ for the equation (\ref{2.16}). Substituting the unique solution $R_\ast(t)$ into the expressions (\ref{2.14})-(\ref{2.15}), one finds that the radially symmetric $T$-periodic positive solution of the problem (\ref{1.2})-(\ref{1.6}) is uniquely 
determined.

Let us introduce the concept of linear stability/instability and asymptotic stability of the periodic solution $(\sigma_\ast(r,t),p_\ast(r,t),R_\ast(t))$. By \emph{linear stability/instability}, we mean: \\
Linearize the system (\ref{1.2})-(\ref{1.6}) at  $(\sigma_\ast(r,t),p_\ast(r,t),R_\ast(t))$ by writing
\begin{equation*}
\begin{split}
  \sigma(r,\theta,\varphi,t)=\sigma_\ast(r,t)+\varepsilon w(r,\theta,\varphi,t), \ \ \ \ p(r,\theta,\varphi,t)=p_\ast(r,t)+\varepsilon q(r,\theta,\varphi,t),
\end{split}
\end{equation*}
$$\partial\Omega(t):\ r=R_\ast(t)+\varepsilon\rho(\theta,\varphi,t),$$
and collect only the $\varepsilon$-order terms. Since the original problem is translation invariant, we say that the periodic solution $(\sigma_\ast(r,t),p_\ast(r,t),R_\ast(t))$ is linearly stable in the sense that for any initial data $\rho|_{t=0}=\rho_0(\theta,\varphi)$ and $w|_{t=0}=w_0(r,\theta,\varphi)$,
\begin{equation*}
  \Big|\rho(\theta,\varphi,t)-\sum_{m=-1}^1a_mY_{1,m}(\theta,\varphi)\Big|\leq Ce^{-\delta t}, \ \ \ \ t>\overline{t}
\end{equation*}
for some constants $a_m$, $\delta>0$ and $\overline{t}>0$. The periodic solution is said to be linearly unstable if it is not linearly stable.

By \emph{asymptotic stability}, we mean the following:\\
Given any initial data
\begin{equation}\label{0.68}
\begin{split}
  \partial\Omega(0):\ r=R_0(\theta,\varphi)=R_\ast(0)+\varepsilon\rho_0(\theta,\varphi),\\ \sigma|_{t=0}=\sigma_0(r,\theta,\varphi)=\sigma_\ast(r,0)+\varepsilon w_0(r,\theta,\varphi),
\end{split}
\end{equation}
if $|\varepsilon|$ is sufficiently small, then there exists a solution to (\ref{1.2})-(\ref{1.6}) for all $t>0$, and
\begin{equation*}
  \partial\Omega(t)\ \ \ \mbox{behaves\ like} \ \ \ \partial B_{R_\ast(t)}(a)=\big\{x:|x-a|=R_\ast(t)\big\}, \ \ \ t>\overline{t}
\end{equation*}
for some center $a$ and $\overline{t}>0$.

Set
\begin{equation}\label{0.1}
  \mu_\ast=\frac{\int_0^T\frac{6}{R^3_\ast(t)}dt}{-\frac{\widetilde{\sigma}}{2}\int_0^T\left[R_\ast(t)\frac{I_3(R_\ast(t))}{I_2(R_\ast(t))}
  -R_\ast(t)\frac{I_0(R_\ast(t))}{I_1(R_\ast(t))}+2\right]dt}>0, \ \ \  (\mbox{two-space\ dimension}),
\end{equation}
or,
\begin{equation}\label{0.16}
  \mu_\ast=\frac{\int_0^T\frac{4}{R^3_\ast(t)}dt}{\int_0^T\frac{\widetilde{\sigma}}{3}R^2_\ast(t)\left[P_1(R_\ast(t))-P_2(R_\ast(t))\right]dt}>0,\ \ \  \ (\mbox{three-space\ dimension}).
\end{equation}
It was recently established in \cite{HZH3} that the unique radially symmetric $T$-periodic positive solution $(\sigma_\ast(r,t),p_\ast(r,t),R_\ast(t))$ is linearly stable under non-radial perturbations for $0<\mu<\mu_\ast$ and linearly unstable for $\mu>\mu_\ast$ in the two-space dimensional case; see \cite{HX} for the similar result in the three-space dimensional case. In this paper, we further study the asymptotic stability for the fully nonlinear problem. Compared with the linearized problem, the $O(\varepsilon^2)$ terms cannot be dropped. On the one hand, we need to find the necessary PDE estimates to estimate the nonlinear error terms. On the other hand, since this system is translation invariant in the coordinate space, the center of the limiting sphere is not known in advance which depends on the perturbation of mode 1. To be specific, the perturbation of mode 1 results in the translation of the origin with its magnitude the same order as the perturbation, as well as the decay behavior in time $t$.  This is a challenge for this type of problems. We employ on a fixed point theorem to find the correct translation of the origin which is technical and carried out by Theorem \ref{th3}.

Throughout the paper, $\mu_\ast$ is defined by (\ref{0.16}). The main result of this paper is stated as follows (see Theorem \ref{th4} for more explicit statement).
\begin{thm}
Let $\mu<\mu_\ast$, then the unique radially symmetric $T$-periodic solution $(\sigma_\ast(r,t),p_\ast(r,t),R_\ast(t))$ of the system (\ref{1.2})-(\ref{1.6}) is asymptotically stable modulus translation, i.e., there exists a new center $\varepsilon a^\ast(\varepsilon)$, where $a^\ast(\varepsilon)$ is a bounded function of $\varepsilon$, such that $\partial\Omega(t)$ behaves like
\begin{equation*}
  \partial B_{R_\ast(t)}(\varepsilon a^\ast(\varepsilon))=\big\{x:|x-\varepsilon a^\ast(\varepsilon)|=R_\ast(t)\big\}, \ \ \ \ \ t>\overline{t}
\end{equation*}
for some $\overline{t}>0$.
\end{thm}

\begin{rem}
Note that our method and result is also applicable to the two-space dimensional case.
\end{rem}
The structure of the rest of this paper is arranged as follows. In the next section, we collect some preliminaries which are needed in the sequel. In Section 3, we transform the nonlinear free boundary problem (\ref{1.2})-(\ref{1.6}) into an initial-boundary value problem defined on a ball, which we then rewrite as a inhomogeneous linear system (\ref{3.150})-(\ref{3.190}) in Section 4.

The method we shall use is a fixed point argument:

\noindent(i) For given inhomogeneous terms $(f^i,b^j)$, solve the inhomogeneous linear system (\ref{3.150})-(\ref{3.190}) by using the spherical harmonic expansion;

\noindent(ii) define new inhomogeneous terms $(\widetilde{f}^i,\widetilde{b}^j)$ (see (\ref{0.17}));

\noindent(iii) establish the existence of a fixed point for the mapping:
\begin{equation*}
  S:(f^i,b^j)\rightarrow(\widetilde{f}^i,\widetilde{b}^j).
\end{equation*}
Furthermore, in Section 4, we establish decay estimates for each mono-mode system obtained from the spherical harmonic expansion. In order to derive decay estimates for mode 1 terms, we need to translate the origin $x=0$ to $x=\varepsilon a^\ast(\varepsilon)$, which is carried out by Theorem \ref{th3}. Section 5 aims at showing the mapping $S$ admits a fixed point, thereby completing the proof of 
the asymptotic stability for the fully nonlinear problem.

\section{Preliminaries}
In this section, we collect some preliminaries which are needed in the sequel.

The function $P_n(r)$ introduced by Friedman and Hu \cite{FH1} is defined by
\begin{equation}\label{2.1}
  P_n(r)=\frac{I_{n+3/2}(r)}{rI_{n+1/2}(r)}, \ \ \ \ n=0,1,2,3,\cdots
\end{equation}
where $I_m(r)$ is the modified Bessel function, $m\geq0$ and $r>0$.
Recall \cite{FH1,HZH4} that
\begin{equation}\label{2.2}
  P_0(r)=\frac{1}{r}\coth r-\frac{1}{r^2},
\end{equation}
\begin{equation}\label{2.18}
  P_n(r)=\frac{1}{r^2P_{n+1}(r)+2n+3},
\end{equation}
\begin{equation}\label{0.8}
  P_n(r)>P_{n+1}(r),
\end{equation}
\begin{equation}\label{0.2}
  \frac{d}{dr}\Big(\frac{I_{n+1/2}(r)}{r^{1/2}}\Big)=\frac{I_{n+3/2}(r)+\frac{n}{r}I_{n+1/2}(r)}{r^{1/2}}.
\end{equation}
\begin{lem}
For $n\geq1$ and $r>0$,
\begin{equation}\label{2.17}
  rP_n(r)<rP_0(r)\leq1.
\end{equation}
\end{lem}
\noindent{\bf Proof.}\ By (\ref{2.2}), we get
\begin{equation*}
  \frac{d}{dr}(rP_0(r))=\frac{d}{dr}\Big(\coth r-\frac{1}{r}\Big)=\frac{(e^{2r}-1)^2-4r^2e^{2r}}{r^2(e^{2r}-1)^2}.
\end{equation*}
Denote the numerator by $F(r)$. We compute
\begin{equation*}
  F(0)=0\ \ \ \ \ \mbox{and}\ \ \ \ \frac{dF(r)}{dr}=4e^{2r}(e^{2r}-1-2r-2r^2)>0,
\end{equation*}
then $rP_0(r)$ is an increasing function, and
\begin{equation*}
  rP_0(r)\leq\lim_{r\rightarrow\infty}rP_0(r)=1.
\end{equation*}
It follows from (\ref{0.8}) that (\ref{2.17}) holds. Therefore, our proof is complete.
\hfill $\Box$

\begin{lem}\label{lemma8}
(\mbox{see} \cite[Lemma 4.1]{HX}) The following relations hold:
\begin{equation}\label{0.12}
  \frac{\partial\sigma_\ast}{\partial r}\Big|_{r=R_\ast(t)}=\phi(t)R_\ast(t)P_0(R_\ast(t)),
\end{equation}
\begin{equation}\label{0.13}
  \frac{\partial p_\ast}{\partial r}\Big|_{r=R_\ast(t)}=-\frac{dR_\ast(t)}{dt},
\end{equation}
\begin{equation}\label{0.14}
  \frac{\partial^2p_\ast}{\partial r^2}\Big|_{r=R_\ast(t)}=-\mu\phi(t)R^2_\ast(t)P_0(R_\ast(t))P_1(R_\ast(t))-\frac{1}{R_\ast(t)}\frac{dR_\ast(t)}{dt}.
\end{equation}
\end{lem}

 As established in \cite{CF},
we have the following local existence theorem:
\begin{thm}\label{th1}
If
\begin{equation}\label{0.57}
  (\sigma_0,R_0)\in C^{1+\gamma}(\overline{\Omega}(0))\times C^{4+\alpha}(\partial \Omega(0))\ \ \ \mbox{and}\ \ \sigma_0=\phi(0) \ \ \ \mbox{on}\ \ \partial\Omega(0)
\end{equation}
for some $\alpha$, $\gamma$ $\in(0,1)$, then there exists a unique solution $(\sigma,p,R)$ of (\ref{1.2})-(\ref{1.6}) for $t\in[0,T]$ with some $T>0$, and
\begin{equation*}
  \sigma\in C^{1+\gamma,(1+\gamma)/2}\Big(\bigcup_{t\in[0,T]}\overline{\Omega}(t)\times\{t\}\Big)\cap C^{2+2\alpha/3,1+\alpha/3}\Big(\bigcup_{t\in[t_0,T]}\overline{\Omega}(t)\times\{t\}\Big) \ \ \ \mbox{for\ any}\ t_0>0,
\end{equation*}
\begin{equation*}
  p\in C^{2+\alpha,\alpha/3}\Big(\bigcup_{t\in[0,T]}\overline{\Omega}(t)\times\{t\}\Big), \ \ \ R\in C^{4+\alpha,1+\alpha/3}.
\end{equation*}
\end{thm}

The following fixed point theorem 
(\hspace{-0.07em}\cite{HZH4}) will be  crucial  
in the proof of Theorem \ref{th3}:
\begin{thm}\label{th2}
Let $(X,\|\cdot\|)$ be a Banach space and let $\overline{B}_K(a_0)$ denote the closed ball in $X$ with center $a_0$ and radius $K$. Let $F$ be a mapping from $\overline{B}_K(a_0)$ into $X$ and
\begin{equation*}
F(x)=F_1(x)+\varepsilon G(x),
\end{equation*}
such that

(i) $F_1'(x)$ and $G'(x)$ are both continuous for $x\in\overline{B}_K(a_0)$,

(ii) $F_1(a_0)=0$ and the operator $F'_1(a_0)$ is invertible.

\noindent Then for small $|\varepsilon|$, the equation $F(x)=0$ admits a unique solution $x$ in $\overline{B}_K(a_0)$.
\end{thm}

\section{Transformation}
In this section, we transform the nonlinear free boundary problem into a nonlinear perturbation in a spherical region.

Let us assume the solution of the system (\ref{1.2})-(\ref{1.6}) is of the form
\begin{equation*}
\begin{split}
    \partial\Omega(t):\ r&=R_\ast(t)+\varepsilon\rho(\theta,\varphi,t),\\
    \sigma(r,\theta,\varphi,t)&=\sigma_\ast(r,t)+\varepsilon w(r,\theta,\varphi,t),\\
    p(r,\theta,\varphi,t)&=p_\ast(r,t)+\varepsilon q(r,\theta,\varphi,t).
\end{split}
\end{equation*}
Recall \cite{HZH1} that
\begin{equation*}
  \vec{n}=\frac{1}{\sqrt{1+|\varepsilon\nabla_\omega\rho|^2/(R_\ast(t)+\varepsilon\rho)^2}}
  \Big(\vec{e}_r-\frac{\varepsilon}{R_\ast(t)+\varepsilon\rho}\frac{\partial\rho}{\partial\theta}\vec{e}_\theta
  -\frac{\varepsilon}{(R_\ast(t)+\varepsilon\rho)\sin\theta}\frac{\partial\rho}{\partial\varphi}\vec{e}_\varphi\Big),
\end{equation*}
and
\begin{equation*}
 \nabla=\vec{e}_r\frac{\partial}{\partial r}+\vec{e}_\theta\frac{1}{r}\frac{\partial}{\partial\theta}+\vec{e}_\varphi\frac{1}{r\sin\theta}\frac{\partial}{\partial\varphi}
 =\vec{e}_r\frac{\partial}{\partial r}+\frac{1}{r}\nabla_\omega,
\end{equation*}
where $\nabla_\omega=\vec{e}_\theta\frac{\partial}{\partial\theta}+\vec{e}_\varphi\frac{1}{\sin\theta}\frac{\partial}{\partial\varphi}$.

The system (\ref{1.2})-(\ref{1.6}) can then be written in terms of $(w,q,\rho)$ as follows:
\begin{equation}\label{3.2}
    \Delta w=w \ \ \ \mbox{in} \ \Omega(t),\ t>0,
\end{equation}
\begin{equation}\label{3.4}
    -\Delta q=\mu w\ \ \ \mbox{in} \ \partial\Omega(t),\ t>0,
\end{equation}
\begin{equation}\label{3.6}
    \frac{\partial\rho}{\partial t}=-\frac{1}{\varepsilon}\Big(\frac{\partial p_\ast}{\partial\vec{n}}+\varepsilon\frac{\partial q}{\partial\vec{n}}\Big)\sqrt{1+\frac{|\varepsilon\nabla_\omega\rho|^2}{(R_\ast(t)+\varepsilon\rho)^2}}-\frac{1}{\varepsilon}\frac{dR_\ast(t)}{dt} \ \ \mbox{on}\ \partial\Omega(t),\ t>0,
\end{equation}
\begin{equation}\label{3.3}
    w=-\frac{1}{\varepsilon}\Big[\sigma_\ast(R_\ast(t)+\varepsilon\rho,t)-\sigma_\ast(R_\ast(t),t)\Big] \ \ \ \mbox{on} \ \partial\Omega(t),\ t>0,
\end{equation}
\begin{equation}\label{3.5}
    q=-\frac{1}{\varepsilon}\big[p_\ast(R_\ast(t)+\varepsilon\rho,t)-\kappa\big] \ \ \ \mbox{on} \ \partial\Omega(t),\ t>0.
\end{equation}
By the Taylor expansion, one of the right-hand side terms of (\ref{3.6}) can be written in the following way,
\begin{equation}\label{3.7}
\begin{split}
  \sqrt{1+\frac{|\varepsilon\nabla_\omega\rho|^2}{(R_\ast(t)+\varepsilon\rho)^2}}\frac{\partial p_\ast}{\partial\vec{n}}\Big|_{r=R_\ast(t)+\varepsilon\rho}&=\sqrt{1+\frac{|\varepsilon\nabla_\omega\rho|^2}{(R_\ast(t)+\varepsilon\rho)^2}}\nabla p_\ast|_{r=R_\ast(t)+\varepsilon\rho}\cdot\vec{n}\\
  &=\frac{\partial p_\ast(R_\ast(t)+\varepsilon\rho,t)}{\partial r}\\
  &=\frac{\partial p_\ast(R_\ast(t),t)}{\partial r}+\frac{\partial^2 p_\ast(R_\ast(t),t)}{\partial r^2}\varepsilon\rho+\varepsilon^2 P_\varepsilon.
\end{split}
\end{equation}
In addition,
\begin{equation}\label{3.07}
\begin{split}
  \sqrt{1+\frac{|\varepsilon\nabla_\omega\rho|^2}{(R_\ast(t)+\varepsilon\rho)^2}}\frac{\partial q}{\partial\vec{n}}\Big|_{r=R_\ast(t)+\varepsilon\rho}&=\sqrt{1+\frac{|\varepsilon\nabla_\omega\rho|^2}{(R_\ast(t)+\varepsilon\rho)^2}}\nabla q|_{r=R_\ast(t)+\varepsilon\rho}\cdot\vec{n}\\
  &=\frac{\partial q}{\partial r}\Big|_{r=R_\ast(t)+\varepsilon\rho}-\frac{\varepsilon}{(R_\ast(t)+\varepsilon\rho)^2}\frac{\partial\rho}{\partial\theta}\frac{\partial q}{\partial\theta}\Big|_{r=R_\ast(t)+\varepsilon\rho}\\
  &\quad-\frac{\varepsilon}{(R_\ast(t)+\varepsilon\rho)^2\sin^2\theta}\frac{\partial\rho}{\partial\varphi}\frac{\partial q}{\partial\varphi}\Big|_{r=R_\ast(t)+\varepsilon\rho}.
\end{split}
\end{equation}
Again by the Taylor expansion, the right-hand side terms of (\ref{3.3}) and (\ref{3.5}) can be expanded into  the following format, respectively,
\begin{equation}\label{3.8}
  \sigma_\ast(R_\ast(t)+\varepsilon\rho,t)-\sigma_\ast(R_\ast(t),t)=\frac{\partial\sigma_\ast}{\partial r}\Big|_{r=R_\ast(t)}\varepsilon\rho+\varepsilon^2S_\varepsilon,
\end{equation}
and
\begin{equation}\label{3.9}
\begin{split}
  p_\ast(R_\ast(t)+\varepsilon\rho,t)-\kappa=&p_\ast(R_\ast(t),t)+\frac{\partial p_\ast}{\partial r}\Big|_{r=R_\ast(t)}\varepsilon\rho+\varepsilon^2K_{1\varepsilon}\\
  &-\frac{1}{R_\ast(t)}+\frac{\varepsilon}{R^2_\ast(t)}\Big(\rho+\frac{1}{2}\Delta_\omega\rho\Big)+\varepsilon^2K_{2\varepsilon}\\
  =&\varepsilon\Big[\frac{1}{R_\ast^2(t)}\Big(\rho+\frac{1}{2}\Delta_\omega\rho\Big)+\frac{\partial p_\ast}{\partial r}\Big|_{r=R_\ast(t)}\rho\Big]+\varepsilon^2K_\varepsilon,
\end{split}
\end{equation}
where we have used the fact, see \cite[Theorem 8.1]{FR2},
\begin{equation*}
  \kappa=\frac{1}{R}-\frac{\varepsilon}{R^2}\Big(\rho+\frac{1}{2}\Delta_\omega\rho\Big)+\varepsilon^2K
\end{equation*}
with
\begin{equation*}
    \Delta_\omega\rho=\frac{1}{\sin\theta}\frac{\partial}{\partial\theta}\left(\sin\theta\frac{\partial\rho}{\partial\theta}\right)
    +\frac{1}{\sin^2\theta}\frac{\partial^2\rho}{\partial\varphi^2}.
\end{equation*}

The Hanzawa transformation is defined by
\begin{equation*}
    r=r'+\chi(R-r')\varepsilon \rho(\theta,\varphi,t), \ \ \ t=t', \ \ \  \ \theta=\theta', \ \ \ \varphi=\varphi',
\end{equation*}
here
\begin{equation*}
   \chi(z)\in C^{\infty}, \ \ \ \ \chi(z)=\left\{\begin{aligned}
&0,  \ \ \mbox{if} \ |z|\geq\frac{3}{4}\delta_0 \\
&1, \ \ \mbox{if}\ |z|<\frac{1}{4}\delta_0
\end{aligned}\right., \ \ \ \
 \left|\frac{d^{k}\chi}{dz^{k}}\right|\leq\frac{C}{\delta_0^{k}},\ \ \ \ (\delta_0\ \mbox{positive\ and\ small}).
\end{equation*}
Under the Hanzawa transformation and (\ref{3.7})-(\ref{3.9}), the system (\ref{3.2})-(\ref{3.5}) is transformed into
\begin{equation}\label{3.10}
    -\Delta' w'+w'=\varepsilon A_\varepsilon w' \ \ \ \mbox{in} \ B_{R_\ast(t)},\ t>0,
\end{equation}
\begin{equation}\label{3.12}
    -\Delta' q'-\mu w'=\varepsilon A_\varepsilon q' \ \ \ \mbox{in} \ B_{R_\ast(t)},\ t>0,
\end{equation}
\begin{equation}\label{3.14}
    \frac{\partial\rho'}{\partial t'}=-\frac{\partial^2p_\ast}{\partial r'^2}\rho'-\frac{\partial q'}{\partial r'}+\varepsilon B^1_\varepsilon \ \ \ \mbox{on} \ \partial B_{R_\ast(t)},\ t>0,
\end{equation}
\begin{equation}\label{3.11}
    w'=-\frac{\partial\sigma_\ast}{\partial r'}\rho'+\varepsilon B^2_\varepsilon  \ \ \ \mbox{on} \ \partial B_{R_\ast(t)},\ t>0,
\end{equation}
\begin{equation}\label{3.13}
    q'=-\frac{1}{R_\ast^2(t)}\Big(\rho'+\frac{1}{2}\Delta_\omega\rho'\Big)-\frac{\partial p_\ast}{\partial r'}\rho'+\varepsilon B^3_\varepsilon \ \ \ \mbox{on} \ \partial B_{R_\ast(t)},\ t>0,
\end{equation}
where $B_{R_\ast(t)}$ is the ball with radius $R_\ast(t)$, $A_\varepsilon$ given in \cite{FH2} is a second order differential operator in $(r',\theta',\varphi')$, and
\begin{equation*}
\begin{split}
  B_\varepsilon^1=&-\frac{1}{\varepsilon^2}\frac{\partial p_\ast(R_\ast(t')+\varepsilon\rho',t')}{\partial r'}+\frac{1}{(R_\ast(t')+\varepsilon\rho')^2}\frac{\partial\rho'}{\partial\theta'}\frac{\partial q'}{\partial\theta'}\\
  &+\frac{1}{(R_\ast(t')+\varepsilon\rho')^2\sin^2\theta'}\frac{\partial\rho'}{\partial\varphi'}\frac{\partial q'}{\partial\varphi'}-\frac{1}{\varepsilon^2}\frac{dR_\ast(t')}{dt'}+\frac{1}{\varepsilon}\frac{\partial^2 p_\ast(R_\ast(t'),t')}{\partial r'^2}\rho',
\end{split}
\end{equation*}
\begin{equation*}
\begin{split}
  B^2_\varepsilon=&-\frac{1}{\varepsilon^2}\big[\sigma_\ast(R_\ast(t')+\varepsilon\rho',t')-\sigma_\ast(R_\ast(t'),t')\big]
  +\frac{1}{\varepsilon}\frac{\partial\sigma_\ast(R_\ast(t'),t')}{\partial r'}\rho',
\end{split}
\end{equation*}
\begin{equation*}
  B_\varepsilon^3=-\frac{1}{\varepsilon^2}\big[p_\ast(R_\ast(t')+\varepsilon\rho',t')-\kappa\big]+\frac{1}{\varepsilon R_\ast^2(t')}\Big(\rho'+\frac{1}{2}\Delta_\omega\rho'\Big)+\frac{1}{\varepsilon}\frac{\partial p_\ast(R_\ast(t'),t')}{\partial r'}\rho'.
\end{equation*}

As pointed out in \cite{FH2} that all terms of $A_\varepsilon$ do not involve any singularity, then it follows from Theorem \ref{th1} that, for $T>1$,
\begin{equation*}
  A_\varepsilon w'\in C^{2\alpha/3,\alpha/3}(B_{R_\ast(t)}\times[0,T]),
\end{equation*}
\begin{equation*}
  A_\varepsilon q'\in C^{\alpha,\alpha/3}(B_{R_\ast(t)}\times[0,T]).
\end{equation*}
Notice that the term $\frac{1}{\varepsilon}$ appearing in $A_\varepsilon$ is cancelled out by the coefficient that accompanies it, so that both $C^{2\alpha/3,\alpha/3}$ norm of $A_\varepsilon w'$ and $C^{\alpha,\alpha/3}$ norm of $A_\varepsilon q'$ are uniformly bounded in $\varepsilon$.

Moreover, although $\sin^2\theta'$ appears in the denominator in the last term of $B_\varepsilon^1$, this incurs no singularities. Indeed, one can simply choose a different coordinate system to deal with this problem, i.e., $B_\varepsilon^1$ is a function defined on the unit sphere $\Sigma=\{x:|x|=1\}$ (rather than in the variable $(\theta,\varphi)$). Then (\ref{2.15}) and Theorem \ref{th1} yield $B_\varepsilon^1\in C^{1+\alpha,\alpha/3}(\partial B_{R_\ast(t)}\times[0,T])$. It follows from (\ref{3.7}) and (\ref{0.13}) that $C^{1+\alpha,\alpha/3}$ norm of $B_\varepsilon^1$ is uniformly bounded in $\varepsilon$.
By (\ref{2.14}) and Theorem \ref{th1}, we find that $B_\varepsilon^2\in C^{1+2\alpha/3,1+\alpha/3}(\partial B_{R_\ast(t)}\times[0,T])$, and (\ref{3.8}) implies that $C^{1+2\alpha/3,1+\alpha/3}$ norm of $B_\varepsilon^2$ is uniformly bounded in $\varepsilon$.
Similarly, by (\ref{3.9}),
we immediately derive that $B_\varepsilon^3\in C^{2+\alpha,\alpha/3}(\partial B_{R_\ast(t)}\times[0,T])$, and its $C^{2+\alpha,\alpha/3}$ norm is uniformly bounded in $\varepsilon$.

For notational convenience, we shall denote functions $w'(r',\theta',\varphi',t')$, $q'(r',\theta',\varphi',t')$ and $\rho'(\theta',\varphi',t')$ again by $w(r,\theta,\varphi,t)$, $q(r,\theta,\varphi,t)$ and $\rho(\theta,\varphi,t)$, respectively, in the rest of this paper.

\section{The Inhomogeneous Linear System}
In this section, we rewrite (\ref{3.10})-(\ref{3.13}) as a inhomogeneous linear system (\ref{3.150})-(\ref{3.190}), and use spherical harmonic expansions to solve this inhomogeneous system. Furthermore, we establish decay estimates for each mono-mode system obtained from the spherical harmonic expansion.

Specifically, we consider the system (\ref{3.10})-(\ref{3.13}) where the $\varepsilon$ terms of any order are replaced by given functions,
%
\begin{equation}\label{3.150}
    -\Delta w+w=\varepsilon f^1(r,\theta,\varphi,t) \ \ \ \mbox{in} \ B_{R_\ast(t)},\ t>0,
\end{equation}
\begin{equation}\label{3.160}
    -\Delta q=\mu w+\varepsilon f^2(r,\theta,\varphi,t) \ \ \ \mbox{in} \ B_{R_\ast(t)},\ t>0,
\end{equation}
\begin{equation}\label{3.170}
\begin{split}
    \frac{\partial\rho}{\partial t}=&-\frac{\partial^2p_\ast}{\partial r^2}\Big|_{r=R_\ast(t)}\rho-\frac{\partial q}{\partial r}\Big|_{r=R_\ast(t)}+\varepsilon b^1(\theta,\varphi,t) \ \ \ \ \mbox{on} \ \partial B_{R_\ast(t)}, \ t>0,
\end{split}
\end{equation}
\begin{equation}\label{3.180}
     w=-\frac{\partial\sigma_\ast}{\partial r}\Big|_{r=R_\ast(t)}\rho+\varepsilon b^2(\theta,\varphi,t)  \ \ \ \ \mbox{on} \ \partial B_{R_\ast(t)}, \ t>0,
\end{equation}
\begin{equation}\label{3.190}
\begin{split}
    q=&-\frac{1}{R_\ast^2(t)}\Big(\rho+\frac{1}{2}\Delta_\omega\rho\Big)-\frac{\partial p_\ast}{\partial r}\Big|_{r=R_\ast(t)}\rho+\varepsilon b^3(\theta,\varphi,t) \ \ \ \mbox{on} \ \partial B_{R_\ast(t)}, \ t>0,
\end{split}
\end{equation}
and we add initial conditions (cf. (\ref{0.68}))
\begin{equation}\label{3.20}
  w|_{t=0}=w_0(r,\theta,\varphi), \ \ \ \quad \ \ \ \rho|_{t=0}=\rho_0(\theta,\varphi).
\end{equation}

Assume that the functions $f^i$, $b^j$ satisfy the following properties
\begin{equation}\label{4.1}
  \sqrt{|\varepsilon|}\left(\int_0^\infty e^{2\delta_1t}\|f^1(\cdot,t)\|^2_{L^2(B_{R_\ast(t)})}dt\right)^{1/2}\leq1,
\end{equation}
\begin{equation}\label{4.2}
  \sqrt{|\varepsilon|}\left(\int_0^\infty e^{2\delta_1t}\|f^2(\cdot,t)\|^2_{L^2(B_{R_\ast(t)})}dt\right)^{1/2}\leq1,
\end{equation}
\begin{equation}\label{4.3}
  \sqrt{|\varepsilon|}\left(\int_0^\infty e^{2\delta_1t}\|b^1(\cdot,t)\|^2_{H^{1/2}(\partial B_{R_\ast(t)})}dt\right)^{1/2}\leq1,
\end{equation}
\begin{equation}\label{4.4}
  \sqrt{|\varepsilon|}\left(\int_0^\infty e^{2\delta_1t}\|b^2(\cdot,t)\|^2_{H^{1}(\partial B_{R_\ast(t)})}dt\right)^{1/2}\leq1,
\end{equation}
\begin{equation}\label{4.5}
  \sqrt{|\varepsilon|}\left(\int_0^\infty e^{2\delta_1t}\|b^3(\cdot,t)\|^2_{H^{3/2}(\partial B_{R_\ast(t)})}dt\right)^{1/2}\leq1,
\end{equation}
where $\delta_1>0$ is sufficiently small, and for some $\alpha\in(0,1)$,
\begin{equation}\label{4.6}
  \sqrt{|\varepsilon|}\|f^1\|_{C^{2\alpha/3,\alpha/3}(B_{R_\ast(t)}\times[0,\infty))}\leq1,
\end{equation}
\begin{equation}\label{4.7}
  \sqrt{|\varepsilon|}\|f^2\|_{C^{\alpha,\alpha/3}(B_{R_\ast(t)}\times[0,\infty))}\leq1,
\end{equation}
\begin{equation}\label{4.8}
  \sqrt{|\varepsilon|}\|b^1\|_{C^{1+\alpha,\alpha/3}(\partial B_{R_\ast(t)}\times[0,\infty))}\leq1,
\end{equation}
\begin{equation}\label{4.9}
  \sqrt{|\varepsilon|}\|b^2\|_{C^{2+2\alpha/3,1+\alpha/3}(\partial B_{R_\ast(t)}\times[0,\infty))}\leq1,
\end{equation}
\begin{equation}\label{4.10}
  \sqrt{|\varepsilon|}\|b^3\|_{C^{2+\alpha,\alpha/3}(\partial B_{R_\ast(t)}\times[0,\infty))}\leq1.
\end{equation}
For simplicity, we also assume that
\begin{equation}\label{4.12}
  \rho_0\in C^{4+\alpha}(\overline{B}_{R_\ast(0)}), \ \ \ \ \ w_0\in C^{2+2\alpha/3}(\overline{B}_{R_\ast(0)}),
\end{equation}
and the compatibility condition of order 2 for $w$ is satisfied.
\begin{rem}\label{rem1}
For given functions $f^i$, $b^j$ 
in some subset $X_1$ of a Banach space $X$, we solve the inhomogeneous linear system (\ref{3.150})-(\ref{3.190}), and derive the estimate of $(w,q,\rho)$. Then we define the new functions $\widetilde{f}^i$, $\widetilde{b}^j$ by
\begin{equation}\label{0.17}
  \widetilde{f}^1=A_\varepsilon w, \ \ \ \ \widetilde{f}^2=A_\varepsilon q,\ \ \ \
  \widetilde{b}^1=B_\varepsilon^1, \ \ \ \widetilde{b}^2=B_\varepsilon^2, \ \ \ \widetilde{b}^3=B_\varepsilon^3.
\end{equation}
We shall show that the mapping $S:(f^i,b^j)\rightarrow(\widetilde{f}^i,\widetilde{b}^j)$ admits a fixed point which is carried out in Section 5, leading to the asymptotic stability.
\end{rem}
\begin{rem}
The initial data are supposed to satisfy the conditions of Theorem \ref{th1} rather than (\ref{4.12}). However, in the process of showing that the mapping $S$ admits a fixed point, we shall need H\"{o}lder estimates of $D_x^2\sigma$ which requires the consistency condition of order 2 for $\sigma$ at $\partial\Omega(0)$. This is quite restrictive. As stated in \cite[Remark 3.2]{FH2}, we can avoid it by taking the initial time at $t=\frac{T}{2}$ instead of at $t=0$. To simplify notation, we shall denote $t=\frac{T}{2}$ and the initial data $(w|_{t=T/2},\rho|_{t=T/2})$ by $t=0$ and $(w_0,\rho_0)$, respectively. Then the consistency condition of order 2 is satisfied at $t=0$, with $\rho_0\in C^{4+\alpha}(\overline{B}_{R_\ast(0)})$ and $w_0\in C^{2+2\alpha/3}(\overline{B}_{R_\ast(0)})$ (cf. (\ref{4.12})).
\end{rem}

Now we proceed to use spherical harmonic expansions to solve the inhomogeneous linear system (\ref{3.150})-(\ref{3.190}). We first formally expand all the functions $f^i$, $b^j$ in terms of spherical harmonics
\begin{equation*}
  f^i(r,\theta,\varphi,t)=\sum^\infty_{n=0}\sum^n_{m=-n}f^i_{n,m}(r,t)Y_{n,m}(\theta,\varphi), \ \ \ \ i=1,2,
\end{equation*}
\begin{equation*}
  b^j(\theta,\varphi,t)=\sum^\infty_{n=0}\sum^n_{m=-n}b^j_{n,m}(t)Y_{n,m}(\theta,\varphi), \ \ \ \ j=1,2,3.
\end{equation*}
Recall \cite{FH2} that (\ref{4.1})-(\ref{4.5}) imply
\begin{equation}\label{4.19}
  |\varepsilon|\int_0^\infty e^{2\delta_1t}\|f_{n,m}^1(\cdot,t)\|^2_{L^2(B_{R_\ast(t)})}dt=F_{n,m}^1, \ \ \ \ \ \sum_{n,m}F_{n,m}^1\leq1,
\end{equation}
\begin{equation}\label{4.20}
  |\varepsilon|\int_0^\infty e^{2\delta_1t}\|f_{n,m}^2(\cdot,t)\|^2_{L^2(B_{R_\ast(t)})}dt=F_{n,m}^2, \ \ \ \ \ \sum_{n,m}F_{n,m}^2\leq1,
\end{equation}
\begin{equation}\label{4.21}
  |\varepsilon|(n+1)\int_0^\infty e^{2\delta_1t}|b_{n,m}^1(t)|^2dt=B_{n,m}^1, \ \ \ \ \ \sum_{n,m}B_{n,m}^1\leq C,
\end{equation}
\begin{equation}\label{4.22}
  |\varepsilon|(n+1)^2\int_0^\infty e^{2\delta_1t}|b_{n,m}^2(t)|^2dt=B_{n,m}^2, \ \ \ \ \ \sum_{n,m}B_{n,m}^2\leq C,
\end{equation}
\begin{equation}\label{4.23}
  |\varepsilon|(n+1)^3\int_0^\infty e^{2\delta_1t}|b_{n,m}^3(t)|^2dt=B_{n,m}^3, \ \ \ \ \ \sum_{n,m}B_{n,m}^3\leq C.
\end{equation}

 We then look for a solution of the following form:
\begin{eqnarray*}
  w(r,\theta,\varphi,t) &=& \sum^\infty_{n=0}\sum^n_{m=-n}w_{n,m}(r,t)Y_{n,m}(\theta,\varphi), \\
   q(r,\theta,\varphi,t) &=& \sum^\infty_{n=0}\sum^n_{m=-n}q_{n,m}(r,t)Y_{n,m}(\theta,\varphi), \\
\rho(\theta,\varphi,t) &=& \sum^\infty_{n=0}\sum^n_{m=-n}\rho_{n,m}(t)Y_{n,m}(\theta,\varphi).
\end{eqnarray*}
By the relation
\begin{equation*}\label{1.20}
    \Delta_\omega Y_{n,m}(\theta,\varphi)+n(n+1)Y_{n,m}(\theta,\varphi)=0\ \ \ \mbox{and}\ \ \ \Delta=\frac{1}{r^2}\frac{\partial}{\partial r}\Big(r^2\frac{\partial}{\partial r}\Big)+\frac{1}{r^2}\Delta_\omega,
\end{equation*}
we obtain that $w_{n,m}(r,t)$, $q_{n,m}(r,t)$ and $\rho_{n,m}(t)$ satisfy
\begin{equation}\label{4.13}
    -\Delta w_{n,m}(r,t)+\left(\frac{n(n+1)}{r^2}+1\right)w_{n,m}(r,t)=\varepsilon f^1_{n,m}(r,t) \ \ \ \mbox{in}\ B_{R_\ast(t)}, \ t>0,
\end{equation}
\begin{equation}\label{4.15}
    -\Delta q_{n,m}(r,t)+\frac{n(n+1)}{r^2}q_{n,m}(r,t)=\mu w_{n,m}(r,t)+\varepsilon f^2_{n,m}(r,t) \ \ \ \mbox{in} \ B_{R_\ast(t)},\ t>0,
\end{equation}
\begin{equation}\label{4.17}
\begin{split}
    \frac{d\rho_{n,m}(t)}{dt}=&-\frac{\partial^2p_\ast}{\partial r^2}\Big|_{r=R_\ast(t)}\rho_{n,m}(t)-\frac{\partial q_{n,m}(r,t)}{\partial r}\Big|_{r=R_\ast(t)}+\varepsilon b^1_{n,m}(t), \ \ \ t>0,
\end{split}
\end{equation}
\begin{equation}\label{4.14}
    w_{n,m}(R_\ast(t),t)=-\frac{\partial\sigma_\ast}{\partial r}\Big|_{r=R_\ast(t)}\rho_{n,m}(t)+\varepsilon b^2_{n,m}(t),  \ \ \  \ \ t>0,
\end{equation}
\begin{equation}\label{4.16}
\begin{split}
    q_{n,m}(R_\ast(t),t)=&-\frac{1}{R_\ast^2(t)}\left(1-\frac{n(n+1)}{2}\right)\rho_{n,m}(t)\\
    &-\frac{\partial p_\ast}{\partial r}\Big|_{r=R_\ast(t)}\rho_{n,m}(t)+\varepsilon b^3_{n,m}(t), \  \ \ t>0,
\end{split}
\end{equation}
\begin{equation}\label{4.18}
    \rho_{n,m}|_{t=0}=\rho_{0,n,m}, \ \ \ \ \ w_{n,m}|_{t=0}=w_{0,n,m}(r) \ \ \ \ \mbox{in}\ \ B_{R_\ast(0)}.
\end{equation}

As in \cite{HZH1}, we can solve (\ref{4.13}) and (\ref{4.14}) in the form
\begin{equation}\label{4.29}
\begin{split}
   w_{n,m}(r,t)=&\Big[-\frac{\partial\sigma_\ast}{\partial r}\Big|_{r=R_\ast(t)}\rho_{n,m}(t)+\varepsilon b^2_{n,m}(t)\Big]\frac{R_\ast^{1/2}(t)}{I_{n+1/2}(R_\ast(t))}\frac{I_{n+1/2}(r)}{r^{1/2}}\\
    &+\varepsilon\xi_{1,n,m}(r,t),
\end{split}
\end{equation}
where $\xi_{1,n,m}(r,t)$ is the solution of
\begin{equation}\label{4.30}
\begin{split}
  -\Delta\xi_{1,n,m}(r,t)+\left(\frac{n(n+1)}{r^2}+1\right)\xi_{1,n,m}(r,t) &= f^1_{n,m}(r,t) \ \ \ \ \mbox{in}\ B_{R_\ast(t)}, \\
  \xi_{1,n,m}(R_\ast(t),t) &= 0. 
\end{split}
\end{equation}

Let
\begin{equation}\label{4.31}
    \psi(r,t)=q_{n,m}(r,t)+\mu w_{n,m}(r,t),
\end{equation}
then $\psi(r,t)$ satisfies
\begin{equation}\label{4.32}
    -\Delta\psi(r,t)+\frac{n(n+1)}{r^2}\psi(r,t)=\varepsilon\mu f^1_{n,m}(r,t)+\varepsilon f^2_{n,m}(r,t) \ \ \ \ \mbox{in}\ B_{R_\ast(t)},
\end{equation}
and by (\ref{4.14}) and (\ref{4.16}), $\psi(r,t)$ satisfies the following boundary condition
\begin{equation}\label{4.33}
\begin{split}
    \psi(R_\ast(t),t)=&\frac{1}{R_\ast^2(t)}\left(\frac{n(n+1)}{2}-1\right)\rho_{n,m}(t)-\frac{\partial p_\ast}{\partial r}\Big|_{r=R_\ast(t)}\rho_{n,m}(t)\\
    &-\mu\frac{\partial\sigma_\ast}{\partial r}\Big|_{r=R_\ast(t)}\rho_{n,m}(t)+\varepsilon b^3_{n,m}(t)+\mu\varepsilon b^2_{n,m}(t).
\end{split}
\end{equation}
The solution of the problem (\ref{4.32})-(\ref{4.33}) is given by
\begin{equation}\label{4.34}
\begin{split}
    \psi(r,t)=&\frac{r^n}{R^n_\ast(t)}\Big\{\frac{1}{R_\ast^2(t)}\left(\frac{n(n+1)}{2}-1\right)\rho_{n,m}(t)-\frac{\partial p_\ast}{\partial r}\Big|_{r=R_\ast(t)}\rho_{n,m}(t)\\
    &-\mu\frac{\partial\sigma_\ast}{\partial r}\Big|_{r=R_\ast(t)}\rho_{n,m}(t)+\varepsilon b^3_{n,m}(t)+\mu\varepsilon b^2_{n,m}(t)\Big\}+\varepsilon\psi_{1,n,m}(r,t),
\end{split}
\end{equation}
where $\psi_{1,n,m}(r,t)$ is the solution of
\begin{equation}\label{4.35}
\begin{split}
  -\Delta\psi_{1,n,m}(r,t)+\frac{n(n+1)}{r^2}\psi_{1,n,m}(r,t) &= \mu f^1_{n,m}(r,t)+f^2_{n,m}(r,t)\ \ \ \ \mbox{in}\ B_{R_\ast(t)}, \\
  \psi_{1,n,m}(R_\ast(t),t) &= 0. 
\end{split}
\end{equation}

It follows from (\ref{4.31}), (\ref{4.34}), (\ref{4.29}) and (\ref{0.2}) that
\begin{equation*}\label{4.37}
\begin{split}
    \frac{\partial q_{n,m}}{\partial r}\Big|_{r=R_\ast(t)}=&\frac{\partial\psi}{\partial r}\Big|_{r=R_\ast(t)}-\mu\frac{\partial w_{n,m}}{\partial r}\Big|_{r=R_\ast(t)}\\
    =&\frac{n}{R_\ast(t)}\Big\{\frac{1}{R_\ast^2(t)}\left(\frac{n(n+1)}{2}-1\right)\rho_{n,m}(t)-\frac{\partial p_\ast}{\partial r}\Big|_{r=R_\ast(t)}\rho_{n,m}(t)\\
    &-\mu\frac{\partial\sigma_\ast}{\partial r}\Big|_{r=R_\ast(t)}\rho_{n,m}(t)+\varepsilon b^3_{n,m}(t)+\mu\varepsilon b^2_{n,m}(t)\Big\}+\varepsilon\frac{\partial\psi_{1,n,m}}{\partial r}\Big|_{r=R_\ast(t)}\\
    &-\mu\Big[-\frac{\partial\sigma_\ast}{\partial r}\Big|_{r=R_\ast(t)}\rho_{n,m}(t)+\varepsilon b^2_{n,m}(t)\Big]\left\{\frac{I_{n+3/2}(R_\ast(t))}{I_{n+1/2}(R_\ast(t))}+\frac{n}{R_\ast(t)}\right\}\\
    &-\varepsilon\mu\frac{\partial\xi_{1,n,m}}{\partial r}\Big|_{r=R_\ast(t)}.
\end{split}
\end{equation*}
After a direct computation, by (\ref{2.1}), (\ref{0.12}) and (\ref{0.13}), we obtain
\begin{equation*}
\begin{split}
    \frac{\partial q_{n,m}}{\partial r}\Big|_{r=R_\ast(t)}
    =&\frac{n}{R_\ast(t)}\Big\{\frac{1}{R_\ast^2(t)}\left(\frac{n(n+1)}{2}-1\right)\rho_{n,m}(t)-\frac{\partial p_\ast}{\partial r}\Big|_{r=R_\ast(t)}\rho_{n,m}(t)\Big\}\\
    &+\mu\frac{\partial\sigma_\ast}{\partial r}\Big|_{r=R_\ast(t)}R_\ast(t)P_n(R_\ast(t))\rho_{n,m}(t)+\frac{n}{R_\ast(t)}\varepsilon b^3_{n,m}(t)\\
    &-\varepsilon\mu R_\ast(t)P_n(R_\ast(t))b^2_{n,m}(t)
    -\varepsilon\mu\frac{\partial\xi_{1,n,m}}{\partial r}\Big|_{r=R_\ast(t)}+\varepsilon\frac{\partial\psi_{1,n,m}}{\partial r}\Big|_{r=R_\ast(t)}\\
    =&\Big\{\frac{n}{R_\ast^3(t)}\left(\frac{n(n+1)}{2}-1\right)+\frac{n}{R_\ast(t)}\frac{dR_\ast(t)}{dt}\\
    &+\mu\phi(t)R^2_\ast(t)P_0(R_\ast(t))P_n(R_\ast(t))\Big\}\rho_{n,m}(t)+\varepsilon\frac{n}{R_\ast(t)}b^3_{n,m}(t)\\
    &-\varepsilon\mu R_\ast(t)P_n(R_\ast(t))b^2_{n,m}(t)
    -\varepsilon\mu\frac{\partial\xi_{1,n,m}}{\partial r}\Big|_{r=R_\ast(t)}+\varepsilon\frac{\partial\psi_{1,n,m}}{\partial r}\Big|_{r=R_\ast(t)}.
\end{split}
\end{equation*}
Inserting the above expression into (\ref{4.17}) and using the relation (\ref{0.14}), we derive
\begin{equation*}\label{3.41}
\begin{split}
   \frac{d\rho_{n,m}(t)}{dt}
    =&\Big\{-\frac{n-1}{R_\ast(t)}\frac{dR_\ast(t)}{dt}-\frac{n}{R_\ast^3(t)}\Big(\frac{n(n+1)}{2}-1\Big)\\
    &+\mu\phi(t)R^2_\ast(t)P_0(R_\ast(t))\big[P_1(R_\ast(t))-P_n(R_\ast(t))\big]\Big\}\rho_{n,m}(t)\\
    &+\varepsilon\mu\frac{\partial\xi_{1,n,m}}{\partial r}\Big|_{r=R_\ast(t)}-\varepsilon\frac{\partial\psi_{1,n,m}}{\partial r}\Big|_{r=R_\ast(t)}\\
    &+\varepsilon b^1_{n,m}(t)+\varepsilon\mu R_\ast(t)P_n(R_\ast(t))b^2_{n,m}(t)-\varepsilon\frac{n}{R_\ast(t)}b^3_{n,m}(t).
\end{split}
\end{equation*}
Then we solve this ODE in the following form
\begin{equation}\label{4.40}
\begin{split}
  \rho_{n,m}(t)
  =\rho_{0,n,m}e^{\int_0^tH_n(\tau)d\tau}+\varepsilon\int_0^tQ_{n,m}(s)e^{\int_s^tH_n(\tau)d\tau}ds,
\end{split}
\end{equation}
where
\begin{equation}\label{4.41}
\begin{split}
  H_n(\tau)
  =&-\frac{n-1}{R_\ast(\tau)}\frac{dR_\ast(\tau)}{d\tau}-\frac{n}{R_\ast^3(\tau)}\Big(\frac{n(n+1)}{2}-1\Big)\\
  &+\mu\phi(\tau)R^2_\ast(\tau)P_0(R_\ast(\tau))\big[P_1(R_\ast(\tau))-P_n(R_\ast(\tau))\big]\\
  =&\Big\{R_\ast(\tau)\big[P_1(R_\ast(\tau))-P_n(R_\ast(\tau))\big]-\frac{n-1}{R_\ast(\tau)}\Big\}\frac{dR_\ast(\tau)}{d\tau}\\
  &-\frac{n}{R_\ast^3(\tau)}\Big(\frac{n(n+1)}{2}-1\Big)
  +\frac{\mu\widetilde{\sigma}}{3}R_\ast^2(\tau)\big[P_1(R_\ast(\tau))-P_n(R_\ast(\tau))\big],
\end{split}
\end{equation}
here we have used the fact, by (\ref{2.16}),
\begin{equation*}
  \mu\phi(\tau)P_0(R_\ast(\tau))=\frac{1}{R_\ast(\tau)}\frac{dR_\ast(\tau)}{d\tau}+\frac{\mu\widetilde{\sigma}}{3},
\end{equation*}
and
\begin{equation}\label{0.3}
\begin{split}
  Q_{n,m}(s)=&\mu\frac{\partial\xi_{1,n,m}}{\partial r}\Big|_{r=R_\ast(s)}-\frac{\partial\psi_{1,n,m}}{\partial r}\Big|_{r=R_\ast(s)}\\
  &+b_{n,m}^1(s)+\mu R_\ast(s)P_n(R_\ast(s))b_{n,m}^2(s)-\frac{n}{R_\ast(s)}b_{n,m}^3(s).
\end{split}
\end{equation}
Note that $\rho_{M,n,m}\triangleq\rho_{0,n,m}e^{\int_0^tH_n(\tau)d\tau}$ is the function $\rho_{n,m}$ corresponding to the case $f_{n,m}^1\equiv f_{n,m}^2\equiv b_{n,m}^1\equiv b_{n,m}^2\equiv b_{n,m}^3\equiv0$ of the system (\ref{4.13})-(\ref{4.16}). The corresponding $w_{n,m}$ and $q_{n,m}$ are denoted by $w_{M,n,m}$ and $q_{M,n,m}$. Let $w_M=\sum_{n=0}^\infty\sum_{m=-n}^nw_{M,n,m}Y_{n,m}$, $q_M=\sum_{n=0}^\infty\sum_{m=-n}^nq_{M,n,m}Y_{n,m}$ and $\rho_M=\sum_{n=0}^\infty\sum_{m=-n}^n\rho_{M,n,m}Y_{n,m}$.

We shall estimate $\rho_{n,m}(t)$ by dividing into three cases: $n=0$, $n=1$ and $n\geq2$.

\subsection{Case 1: $n=0$}

It follows from (\ref{4.40}) and (\ref{0.3}) that
\begin{equation}\label{4.43}
\begin{split}
  \rho_{0,0}(t)=&\rho_{0,0,0}e^{\int_0^tH_0(\tau)d\tau}+\varepsilon\int_0^tQ_{0,0}(s)e^{\int_s^tH_0(\tau)d\tau}ds\\
  =&\rho_{0,0,0}e^{\int_0^tH_0(\tau)d\tau}+\varepsilon\int_0^t\Big[\mu\frac{\partial\xi_{1,0,0}}{\partial r}\Big|_{r=R_\ast(s)}-\frac{\partial\psi_{1,0,0}}{\partial r}\Big|_{r=R_\ast(s)}\Big]e^{\int_s^tH_0(\tau)d\tau}ds\\
  &+\varepsilon\int_0^tb_{0,0}^1(s)e^{\int_s^tH_0(\tau)d\tau}ds+\varepsilon\int_0^t\mu R_\ast(s)P_0(R_\ast(s))b_{0,0}^2(s)e^{\int_s^tH_0(\tau)d\tau}ds\\
  \equiv& \rho_{M,0,0}+\varepsilon L_1+\varepsilon L_2+\varepsilon L_3,
\end{split}
\end{equation}
where, by (\ref{4.41}),
\begin{equation}\label{4.44}
\begin{split}
  H_0(\tau)=&\Big\{R_\ast(\tau)\big[P_1(R_\ast(\tau))-P_0(R_\ast(\tau))\big]+\frac{1}{R_\ast(\tau)}\Big\}\frac{dR_\ast(\tau)}{d\tau}\\
  &+\frac{\mu\widetilde{\sigma}}{3}R_\ast^2(\tau)\big[P_1(R_\ast(\tau))-P_0(R_\ast(\tau))\big].
\end{split}
\end{equation}

We now proceed to estimate the last three terms on the right-hand of (\ref{4.43}), respectively. Before doing that, we establish
the following lemma.
\begin{lem}\label{lemma1}
For $\mu>0$, there exists a positive number $\delta$, depending on $\mu$ and $R_\ast$, such that the following is true:
\begin{equation}\label{4.45}
  e^{\int_s^tH_0(\tau)d\tau}\leq Ce^{-\delta(t-s)}.
\end{equation}
\end{lem}
\noindent{\bf Proof.}\ For any $t>0$ and $0<s<t$, there exist $m\in \mathbb{N}$ and $z\in[0,T)$ such that $t-s=mT+z$, then
\begin{equation*}
\begin{split}
  &e^{\int_s^t\left\{R_\ast(\tau)\big[P_1(R_\ast(\tau))-P_0(R_\ast(\tau))\big]+\frac{1}{R_\ast(\tau)}\right\}\frac{dR_\ast(\tau)}{d\tau}d\tau}\\
  &\quad=e^{\int_s^{s+z}\left\{R_\ast(\tau)\big[P_1(R_\ast(\tau))-P_0(R_\ast(\tau))\big]+\frac{1}{R_\ast(\tau)}\right\}dR_\ast(\tau)}\leq C.
\end{split}
\end{equation*}
Furthermore, by the boundedness of $R_\ast(t)$ and (\ref{0.8}), we obtain
\begin{equation*}
\begin{split}
  e^{\int_s^tH_0(\tau)d\tau}&\leq Ce^{\int_s^t\frac{\mu\widetilde{\sigma}}{3}R_\ast^2(\tau)\big[P_1(R_\ast(\tau))-P_0(R_\ast(\tau))\big]d\tau}\\
  &\leq Ce^{-\delta(t-s)}
\end{split}
\end{equation*}
for some $\delta>0$. Hence, the proof is complete. \hfill $\Box$

\begin{lem}\label{lemma2}
If
\begin{equation*}
  \int_0^\infty e^{2\delta_1t}|b(t)|^2dt\leq A, \ \ \ \ \ 0<\delta_1<\delta,
\end{equation*}
then
\begin{equation}\label{4.46a}
  \int_0^\infty e^{2\delta_1t}\Big|\int_0^tb(s)e^{\int_s^tH_0(\tau)d\tau}ds\Big|^2dt\leq CA.
\end{equation}
\end{lem}
\noindent{\bf Proof.}\ Lemma \ref{lemma1} yields
\begin{equation*}
\begin{split}
  J_0&\equiv \int_0^\infty e^{2\delta_1t}\Big|\int_0^tb(s)e^{\int_s^tH_0(\tau)d\tau}ds\Big|^2dt\\
  &\leq C\int_0^\infty e^{2\delta_1t}\Big|\int_0^t|b(s)|e^{-\delta(t-s)}ds\Big|^2dt\\
  &=C\int_0^\infty e^{2\delta_1t}\Big|\int_0^t|b(t-s)|e^{-\delta s}ds\Big|^2dt\\
  &\leq C\int_0^\infty\left\{\int_0^t|b(t-s)|^2e^{2\delta_1(t-s)}e^{-\delta s+\delta_1s}ds\int_0^te^{-\delta s+\delta_1s}ds\right\}dt\\
  &\leq C\int_0^\infty\int_0^t|b(t-s)|^2e^{2\delta_1(t-s)}e^{-(\delta-\delta_1)s}dsdt.
\end{split}
\end{equation*}
Changing the order of the integration, we obtain
\begin{equation*}
\begin{split}
  J_0&\leq C\int_0^\infty\int_s^\infty |b(t-s)|^2e^{2\delta_1(t-s)}e^{-(\delta-\delta_1)s}dtds\\
  &=C \int_0^\infty\int_0^\infty |b(t)|^2e^{2\delta_1t}e^{-(\delta-\delta_1)s}dtds\\
  &\leq CA\int_0^\infty e^{-(\delta-\delta_1)s}ds\\
  &\leq CA,
\end{split}
\end{equation*}
which completes our proof.  \hfill $\Box$

\begin{lem}\label{lemma3}
For $0<\delta_1<\delta$, the following estimate holds:
\begin{equation}\label{4.47}
  \int_0^\infty e^{2\delta_1t}\big|L_1\big|^2dt\leq C|\varepsilon|^{-1}\big(F_{0,0}^1+F_{0,0}^2\big).
\end{equation}
\end{lem}
\noindent{\bf Proof.}\ Recall \cite[Lemmas 3.2 and 3.3]{FH1} that
\begin{equation}\label{0.4}
  \left|\frac{\partial\xi_{1,n,m}}{\partial r}\Big|_{r=R_\ast(t)}\right|^2\leq \frac{C}{n+1}\|f_{n,m}^1(\cdot,t)\|^2_{L^2(B_{R_\ast(t)})},
\end{equation}
\begin{equation}\label{0.5}
  \left|\frac{\partial\psi_{1,n,m}}{\partial r}\Big|_{r=R_\ast(t)}\right|^2\leq \frac{C}{n+1}\big(\|f_{n,m}^1(\cdot,t)\|^2_{L^2(B_{R_\ast(t)})}+\|f_{n,m}^2(\cdot,t)\|^2_{L^2(B_{R_\ast(t)})}\big),
\end{equation}
so that, by (\ref{4.19}) and (\ref{4.20}), we obtain
\begin{equation*}
\begin{split}
  &\int_0^\infty e^{2\delta_1t}\Big|\mu\frac{\partial\xi_{1,0,0}}{\partial r}\Big|_{r=R_\ast(t)}-\frac{\partial\psi_{1,0,0}}{\partial r}\Big|_{r=R_\ast(t)}\Big|^2dt\\
  &\quad \leq C\int_0^\infty e^{2\delta_1t}\Big|\frac{\partial\xi_{1,0,0}}{\partial r}\Big|_{r=R_\ast(t)}\Big|^2dt+C\int_0^\infty e^{2\delta_1t}\Big|\frac{\partial\psi_{1,0,0}}{\partial r}\Big|_{r=R_\ast(t)}\Big|^2dt\\
  &\quad \leq C\int_0^\infty e^{2\delta_1t}\|f_{0,0}^1\|^2_{L^2(B_{R_\ast(t)})}dt+C\int_0^\infty e^{2\delta_1t}\big(\|f_{0,0}^1\|^2_{L^2(B_{R_\ast(t)})}+\|f_{0,0}^2\|^2_{L^2(B_{R_\ast(t)})}\big)dt\\
  &\quad \leq C|\varepsilon|^{-1}\big(F_{0,0}^1+F_{0,0}^2\big).
\end{split}
\end{equation*}
Hence, it follows from Lemma \ref{lemma2} that (\ref{4.47}) holds.  \hfill $\Box$

\begin{lem}\label{lemma4}
The following estimates hold:
\begin{equation}\label{4.471}
  \int_0^\infty e^{2\delta_1t}\big|L_2\big|^2dt\leq C|\varepsilon|^{-1}B_{0,0}^1,
\end{equation}
\begin{equation}\label{4.48}
  \int_0^\infty e^{2\delta_1t}\big|L_3\big|^2dt\leq C|\varepsilon|^{-1}B_{0,0}^2.
\end{equation}
\end{lem}
\noindent{\bf Proof.}\  (\ref{4.21}) and Lemma \ref{lemma2} yield (\ref{4.471}). Moreover, by (\ref{2.17}), (\ref{4.22}) and Lemma \ref{lemma2}, we get (\ref{4.48}). Therefore, the proof is complete. \hfill $\Box$

Combining Lemmas \ref{lemma3} and \ref{lemma4}, we establish an estimate for $\rho_{0,0}(t)$ as follows:
\begin{equation}\label{0.10}
  \int_0^\infty e^{2\delta_1t}\big|\rho_{0,0}(t)-\rho_{M,0,0}(t)\big|^2dt\leq C|\varepsilon|\big(F_{0,0}^1+F_{0,0}^2+B_{0,0}^1+B_{0,0}^2\big).
\end{equation}

\subsection{Case 2: $n=1$}

It can be easily seen from (\ref{4.41}) that
\begin{equation*}
\begin{split}
  H_1(\tau)=0.
\end{split}
\end{equation*}
Therefore, (\ref{4.40}) implies that
\begin{equation}\label{3.33}
  \rho_{1,m}(t)=\rho_{0,1,m}+\varepsilon\int_0^tQ_{1,m}(s)ds,
\end{equation}
where
\begin{equation}\label{3.34}
\begin{split}
  Q_{1,m}(s)=&\mu\frac{\partial\xi_{1,1,m}}{\partial r}\Big|_{r=R_\ast(s)}-\frac{\partial\psi_{1,1,m}}{\partial r}\Big|_{r=R_\ast(s)}\\
  &+b_{1,m}^1(s)+\mu R_\ast(s)P_1(R_\ast(s))b_{1,m}^2(s)-\frac{1}{R_\ast(s)}b_{1,m}^3(s).
\end{split}
\end{equation}

As pointed out in the introduction, after the initial values perturbed, say (\ref{0.68}), the domain of the system (\ref{1.2})-(\ref{1.6}) (or (\ref{3.150})-(\ref{3.190})) undergoes a translation owing to the perturbations of mode 1, that is to say, the perturbations of mode 1 terms result in the decay estimates in time $t$, as well as the translation of the origin. In order to establish the asymptotic stability, we need to take care of the center of the limiting sphere. We shall employ on an iterative procedure in which at each iteration we shift the domain to a ``nearly" optimal location with a new center $\varepsilon a=\varepsilon(a_1,a_2,a_3)$. Note that a translation of the origin dose not change the equations (\ref{1.2})-(\ref{1.6}) (or (\ref{3.150})-(\ref{3.190})), but change the initial data. As stated in \cite[Section 6]{FH2}, in the new coordinate system, the initial data are changed as follows, $m=-1,0,1$,
\begin{equation}\label{6.2}
\begin{split}
    w_{0,1,m}\ \mbox{is\ replaced\ by}\ w_{0,1,m}+b_{m+2}\frac{\partial\sigma_\ast(r,0)}{\partial r}+\varepsilon A_m,\\
     \mbox{and}\ \rho_{0,1,m}\ \mbox{is\ replaced\ by}\ \rho_{0,1,m}-b_{m+2}+\varepsilon B_m,
\end{split}
\end{equation}
where $A_m$ and $B_m$ are bounded functions of $(a,\varepsilon)$, and $b_{m+2}$ satisfies
\begin{equation*}
 b_1-b_3=a_1\sqrt{\frac{8\pi}{3}}, \ \ \ \ i(b_1+b_3)=-a_2\sqrt{\frac{8\pi}{3}}, \ \ \ \ b_2=a_3\sqrt{\frac{4\pi}{3}}.
\end{equation*}
Then, in the new coordinate system, the expression $\rho_{1,m}(t)$ is given by
\begin{equation}\label{3.35}
  \rho_{1,m}(t)=\rho_{0,1,m}-b_{m+2}+\varepsilon B_m+\varepsilon\int_0^tQ_{1,m}(s)ds.
\end{equation}

We note that as we go to the system (\ref{3.150})-(\ref{3.190}), we shall always use the translated initial data. As in \cite[Page 632]{FH2}, we take $f^i$, $b^j$ in (\ref{3.150})-(\ref{3.190}) to be functions not only of $(r,\theta,\varphi,t)$, but also of the center $a$, which aims to accommodate the consistency condition of order 2. The function $Q_{1,m}(s)$ depends on $a$ implicitly through the dependence of $\xi_{1,1,m}$, $\psi_{1,1,m}$ and $b^j_{1,m}$ on $a$.

In order to find the center of the limiting sphere and establish the decay estimates in time $t$, we rewrite $\rho_{1,m}(t)$ as
\begin{equation}\label{3.36}
\begin{split}
  \rho_{1,m}(t)&=\rho_{0,1,m}-b_{m+2}+\varepsilon B_m+\varepsilon\int_0^\infty Q_{1,m}(s)ds-\varepsilon\int_t^\infty Q_{1,m}(s)ds\\
  &\triangleq F_m(a)-\varepsilon\int_t^\infty Q_{1,m}(s)ds.
\end{split}
\end{equation}

Let $F(a)=(F_{-1}(a),F_0(a),F_1(a))$.
We have the following critical theorem.
\begin{thm}\label{th3}
For $|\varepsilon|$ small, there exists a new center $\varepsilon a^\ast(\varepsilon)$ such that
\begin{equation*}
  F(a^\ast(\varepsilon))=0.
\end{equation*}
\end{thm}
\noindent{\bf Proof.}\ By (\ref{3.36}), we have
\begin{equation*}
\begin{split}
  F(a)
  = E(a)+\varepsilon G(a).
\end{split}
\end{equation*}
Clearly, there exists $a_0$ such that $E(a_0)=0$ and $E'(a_0)$ is invertible. Moreover, $E'(a)$ and $G'(a)$ are continuous for $a\in \overline{B}_1(a_0)$. Hence, by  applying  Theorem \ref{th2}, we complete the proof.
\hfill $\Box$

Now we proceed to estimate the last term of the right-hand side of (\ref{3.36}). Recall  \cite[Lemma 5.4]{FH2} that
\begin{equation*}
  \int_0^\infty e^{2\delta_1t}\Big|\int_t^\infty b(\tau)d\tau\Big|^2dt\leq\frac{1}{\delta_1^2}\int_0^\infty e^{2\delta_1t}|b(t)|^2dt,
\end{equation*}
so that, by (\ref{2.17}), (\ref{0.4}), (\ref{0.5}) and (\ref{4.19})-(\ref{4.23}), 
we  have
the following estimate:
\begin{equation*}\label{3.37}
\begin{split}
  &\int_0^\infty e^{2\delta_1t}\Big|\int_t^\infty Q_{1,m}(s)ds\Big|^2dt\\
  &\leq C\int_0^\infty e^{2\delta_1t}\big|Q_{1,m}(t)\big|^2dt\\
  &\leq C\int_0^\infty e^{2\delta_1t}\Big(\Big|\frac{\partial\xi_{1,1,m}}{\partial r}\Big|_{r=R_\ast(t)}\Big|^2+\Big|\frac{\partial\psi_{1,1,m}}{\partial r}\Big|_{r=R_\ast(t)}\Big|^2+|b_{1,m}^1|^2+|b_{1,m}^2|^2+|b_{1,m}^3|^2\Big)dt\\
  &\leq C\int_0^\infty e^{2\delta_1t}\Big(\|f_{1,m}^1\|^2_{L^2}+\|f_{1,m}^2\|^2_{L^2}+|b_{1,m}^1|^2+|b_{1,m}^2|^2+|b_{1,m}^3|^2\Big)dt\\
  &\leq C|\varepsilon|^{-1}\big(F_{1,m}^1+F_{1,m}^2+B_{1,m}^1+B_{1,m}^2+B_{1,m}^3\big).
\end{split}
\end{equation*}
It follows that
\begin{equation}\label{0.9}
\begin{split}
  \int_0^\infty e^{2\delta_1t}|\rho_{1,m}(t)-F_m(a^\ast(\varepsilon))|^2&=|\varepsilon|^2\int_0^\infty e^{2\delta_1t}\Big|\int_t^\infty Q_{1,m}(s)ds\Big|^2dt\\
  &\leq C|\varepsilon|\big(F_{1,m}^1+F_{1,m}^2+B_{1,m}^1+B_{1,m}^2+B_{1,m}^3\big).
\end{split}
\end{equation}

\subsection{Case 3: $n\geq2$}

In this subsection, we shall deal with the case $n\geq2$. We first establish
the following lemma:
\begin{lem}\label{lemma5}
Let $\mu<\mu_\ast$, where $\mu_\ast$ is defined by (\ref{0.16}). For $n\geq2$, there exists a small positive number $\delta$, independent of $n$, such that
\begin{equation}\label{4.49}
  e^{\int_s^tH_n(\tau)d\tau}\leq e^{-\delta(n^3+1)(t-s)}.
\end{equation}
\end{lem}
\noindent{\bf Proof.}\ 
 Since $0<\mu<\mu_\ast$ and $R_\ast(t)$ is periodic and bounded, 
it follows from (\ref{4.41}) and (\ref{2.17}) that there exists $\gamma>0$ ($\gamma$ is independent of $n$) such that
\begin{equation*}
\begin{split}
  H_n(\tau)
  \leq-\gamma(n^3+1)
\end{split}
\end{equation*}
for $n$ sufficiently large, say $n>n_0$, which yields (\ref{4.49}).

Recall \cite[Lemma 4.4]{HX} that
\begin{equation*}
  \mu_\ast\leq\frac{\int_0^{T}\frac{j}{R_\ast^3(\tau)}\Big(\frac{j(j+1)}{2}-1\Big)d\tau}{\int_{0}^{T}\frac{\widetilde{\sigma}}{3}
  R_\ast^2(\tau)\big[P_1(R_\ast(\tau))-P_j(R_\ast(\tau))\big]d\tau}, \ \ \ \ \ j\geq2.
\end{equation*}
Then for each fixed $j\in[2,n_0]$ and $0<\mu<\mu_\ast$, we have
\begin{equation}\label{3.42}
\begin{split}
  \int_{s}^{s+T}&\frac{j}{R_\ast^3(\tau)}\Big(\frac{j(j+1)}{2}-1\Big)d\tau-\mu\int_{s}^{s+T}\frac{\widetilde{\sigma}}{3}
  R_\ast^2(\tau)\big[P_1(R_\ast(\tau))-P_j(R_\ast(\tau))\big]d\tau>\gamma_j,
\end{split}
\end{equation}
where $\gamma_j>0$ is small enough, and dependent of $(R_\ast,T,\mu)$ and $j$.
For $t>s$, there exists $m\in\mathbb{N}$ and $\nu\in[0,T)$ such that $t-s=mT+\nu$. Furthermore, by (\ref{4.41}) and (\ref{3.42}), we obtain
\begin{equation*}
\begin{split}
  e^{\int_s^tH_n(\tau)d\tau}=&e^{\int_s^{s+\nu}H_n(\tau)d\tau+m\int_s^{s+T}H_n(\tau)d\tau}\\
  =&\mbox{exp}\Big\{\int_s^{s+\nu}\Big[R_\ast(\tau)\big[P_1(R_\ast(\tau))-P_j(R_\ast(\tau))\big]-\frac{j-1}{R_\ast(\tau)}\Big]\frac{dR_\ast(\tau)}{d\tau}d\tau\\
  &-\int_s^{s+\nu}\frac{j}{R_\ast^3(\tau)}\Big(\frac{j(j+1)}{2}-1\Big)d\tau\\
  &+\mu\int_s^{s+\nu}\frac{\widetilde{\sigma}}{3} R_\ast^2(\tau)\big[P_1(R_\ast(\tau))-P_j(R_\ast(\tau))\big]d\tau\\
  &-m\int_s^{s+T}\frac{j}{R_\ast^3(\tau)}\Big(\frac{j(j+1)}{2}-1\Big)d\tau\\
  &+m\mu\int_s^{s+T}\frac{\widetilde{\sigma}}{3}R_\ast^2(\tau)\big[P_1(R_\ast(\tau))-P_j(R_\ast(\tau))\big]d\tau\Big\}\\
   \leq& e^{-m\widetilde{\delta}_j}=e^{-\frac{t-s-\nu}{T}\widetilde{\delta}_j}\leq e^{-\delta_j(t-s)}
\end{split}
\end{equation*}
for some $\widetilde{\delta}_j>0$ and $\delta_j>0$ (both $\widetilde{\delta}_j$ and $\delta_j>0$ depend on $(R_\ast,T,\mu)$ and $j$).

Let $\delta=\min\{\gamma,\frac{\delta_2}{2^3+1},\frac{\delta_3}{3^3+1},\cdots,\frac{\delta_{n_0}}{n_0^3+1}\}$, then for $n\geq2$, we get (\ref{4.49}),  and this completes the proof.
\hfill $\Box$

\begin{lem}\label{lemma6}
If
\begin{equation*}
  \int_0^\infty e^{2\delta_1t}|b(t)|^2dt\leq A, \ \ \ \ \ 0<\delta_1<\delta,
\end{equation*}
then, for $n\geq2$,
\begin{equation}\label{4.46}
  \int_0^\infty e^{2\delta_1t}\Big|\int_0^tb(s)e^{\int_s^tH_n(\tau)d\tau}ds\Big|^2dt\leq CA(n+1)^{-6},
\end{equation}
where $C$ is a constant independent of $n$.
\end{lem}
\noindent{\bf Proof.}\ It follows from Lemma \ref{lemma5} that
\begin{equation*}
\begin{split}
  J_1&\equiv \int_0^\infty e^{2\delta_1t}\Big|\int_0^tb(s)e^{\int_s^tH_n(\tau)d\tau}ds\Big|^2dt\\
  &\leq \int_0^\infty e^{2\delta_1t}\Big(\int_0^t|b(s)|e^{-\delta(n^3+1)(t-s)}ds\Big)^2dt\\
  &=\int_0^\infty e^{2\delta_1t}\Big(\int_0^t|b(t-\tau)|e^{-\delta(n^3+1)\tau}d\tau\Big)^2dt\\
  &\leq \int_0^\infty\left\{\int_0^t|b(t-\tau)|^2e^{2\delta_1(t-\tau)}e^{-\delta(n^3+1)\tau+\delta_1\tau}d\tau\int_0^te^{-\delta (n^3+1)\tau+\delta_1\tau}d\tau\right\}dt\\
  &\leq C(n+1)^{-3}\int_0^\infty\int_0^t|b(t-\tau)|^2e^{2\delta_1(t-\tau)}e^{-\delta (n^3+1)\tau+\delta_1\tau}d\tau dt.
\end{split}
\end{equation*}
Changing the order of the integration, we obtain
\begin{equation*}
\begin{split}
  J_1&\leq C(n+1)^{-3}\int_0^\infty\int_\tau^\infty |b(t-\tau)|^2e^{2\delta_1(t-\tau)}e^{-\delta (n^3+1)\tau+\delta_1\tau}dtd\tau \\
  &=C(n+1)^{-3} \int_0^\infty\int_0^\infty |b(t)|^2e^{2\delta_1t}e^{-\delta (n^3+1)\tau+\delta_1\tau}dtd\tau \\
  &\leq CA(n+1)^{-3}\int_0^\infty e^{-\delta (n^3+1)\tau+\delta_1\tau}d\tau \\
  &\leq CA(n+1)^{-6},
\end{split}
\end{equation*}
which completes our proof.  \hfill $\Box$

\begin{lem}\label{lemma7}
For all $n\geq2$ and $|m|\leq n$,
\begin{equation*}\label{4.470}
  \int_0^\infty e^{2\delta_1t}\Big|\rho_{n,m}-\rho_{M,n,m}\Big|^2dt\leq C|\varepsilon|(n+1)^{-7}\big(F_{n,m}^1+F_{n,m}^2+B_{n,m}^1+B_{n,m}^2+B_{n,m}^3\big),
\end{equation*}
where $C$ is a constant independent of $n$.
\end{lem}
\noindent{\bf Proof.}\ By (\ref{0.3}), (\ref{2.17}), (\ref{0.4}), (\ref{0.5}) and (\ref{4.19})-(\ref{4.23}), we get
\begin{equation*}
\begin{split}
  &\int_0^\infty e^{2\delta_1t}|Q_{n,m}(t)|^2dt\\
  &\leq C\int_0^\infty e^{2\delta_1t}\Big(\Big|\frac{\partial\xi_{1,n,m}}{\partial r}\Big|_{r=R_\ast(t)}\Big|^2+\Big|\frac{\partial\psi_{1,n,m}}{\partial r}\Big|_{r=R_\ast(t)}\Big|^2+|b_{n,m}^1|^2+|b_{n,m}^2|^2+n^2|b_{n,m}^3|^2\Big)dt\\
  &\leq C\int_0^\infty e^{2\delta_1t}\Big(\frac{1}{n+1}\|f_{n,m}^1\|^2_{L^2}+\frac{1}{n+1}\|f_{n,m}^2\|^2_{L^2}+|b_{n,m}^1|^2+|b_{n,m}^2|^2+n^2|b_{n,m}^3|^2\Big)dt\\
  &\leq C|\varepsilon|^{-1}(n+1)^{-1}\big(F_{n,m}^1+F_{n,m}^2+B_{n,m}^1+B_{n,m}^2+B_{n,m}^3\big).
\end{split}
\end{equation*}
Furthermore, it follows from Lemma \ref{lemma6} that
\begin{equation*}
\begin{split}
  \int_0^\infty e^{2\delta_1t}\Big|\rho_{n,m}-\rho_{M,n,m}\Big|^2dt
  =&|\varepsilon|^2\int_0^\infty e^{2\delta_1t}\Big|\int_0^t Q_{n,m}(s)e^{\int_s^tH_n(\tau)d\tau}ds\Big|^2dt\\
  \leq&C|\varepsilon|(n+1)^{-7}\big(F_{n,m}^1+F_{n,m}^2+B_{n,m}^1+B_{n,m}^2+B_{n,m}^3\big).
\end{split}
\end{equation*}
Hence, the proof is complete.
 \hfill $\Box$

Combining (\ref{0.10}), (\ref{0.9}) and Lemma \ref{lemma7}, we obtain, in the new coordinate system,
\begin{equation}\label{0.11}
  \int_0^\infty e^{2\delta_1t}\|\rho-\rho_M\|^2_{H^{7/2}(\partial B_{R_\ast(t)})}dt\leq C|\varepsilon|.
\end{equation}
\begin{rem}
In fact, $(w_M,q_M,\rho_M)$ is the solution of the corresponding linear system of the original problem (\ref{1.2})-(\ref{1.6}), and its decay estimate in time $t$ has been derived in \cite{HZH3} for the two-space dimensional case and \cite{HX} for the three-space dimensional case, respectively.
\end{rem}

Moreover, we have the following estimates for $w$, $q$ and $\frac{\partial\rho}{\partial t}$.
\begin{lem}\label{th5}
If $\delta_1$ is sufficiently small then
\begin{equation}\label{0.18}
  \int_0^\infty e^{2\delta_1t}\|w-w_M\|^2_{H^{1}( B_{R_\ast(t)})}dt\leq C|\varepsilon|,
\end{equation}
\begin{equation}\label{0.19}
  \int_0^\infty e^{2\delta_1t}\|q-q_M\|^2_{H^{2}( B_{R_\ast(t)})}dt\leq C|\varepsilon|,
\end{equation}
\begin{equation}\label{0.20}
  \int_0^\infty e^{2\delta_1t}\Big\|\frac{\partial}{\partial t}(\rho-\rho_M)\Big\|^2_{H^{1/2}(\partial B_{R_\ast(t)})}dt\leq C|\varepsilon|.
\end{equation}
\end{lem}
\noindent{\bf Proof.}\ Observe that $\widetilde{w}\equiv w_{n,m}-w_{M,n,m}$ satisfies (\ref{4.13}) and (\ref{4.14}) with $\rho_{n,m}$ replaced by $\widetilde{\rho}\equiv\rho_{n,m}-\rho_{M,n,m}$, i.e.,
\begin{equation*}
\begin{split}
  -\Delta \widetilde{w}+\left(\frac{n(n+1)}{r^2}+1\right)\widetilde{w}&=\varepsilon f_{n,m}^1(r,t) \qquad \qquad \quad \ \   \ \ \ \mbox{in} \ B_{R_\ast(t)},\ t>0,\\
  \widetilde{w}(R_\ast(t),t)&=-\lambda\widetilde{\rho}(t)+\varepsilon b_{n,m}^2(t),\  \qquad \ \ \ t>0,
\end{split}
\end{equation*}
where $\lambda=\frac{\partial\sigma_\ast}{\partial r}\Big|_{r=R_\ast(t)}$. 

Define $V=\widetilde{w}-\widetilde{w}|_{r=R_\ast(t)}$, we obtain
\begin{equation*}
\begin{split}
  -\Delta V+\left(\frac{n(n+1)}{r^2}+1\right)V&=\varepsilon f_{n,m}^1+\left(\frac{n(n+1)}{r^2}+1\right)(\lambda\widetilde{\rho}-\varepsilon b_{n,m}^2) \ \ \ \mbox{in} \ B_{R_\ast(t)},\ t>0,\\
  V(R_\ast(t),t)&=0, \qquad \qquad \quad \ t>0.
\end{split}
\end{equation*}
Multiplying both sides with $V$ and integrating over $B_{R_\ast(t)}$, we get
\begin{equation*}
\begin{split}
  &\int_{B_{R_\ast(t)}}|\nabla V|^2dx+\int_{B_{R_\ast(t)}}\left(\frac{n(n+1)}{r^2}+1\right)|V|^2dx\\
  &=\int_{B_{R_\ast(t)}}\varepsilon f_{n,m}^1Vdx+\int_{B_{R_\ast(t)}}\left(\frac{n(n+1)}{r^2}+1\right)\lambda\widetilde{\rho} Vdx-\int_{B_{R_\ast(t)}}\left(\frac{n(n+1)}{r^2}+1\right)\varepsilon b_{n,m}^2Vdx.
\end{split}
\end{equation*}
By Young's inequality, we further have
\begin{equation*}
\begin{split}
  &\int_{B_{R_\ast(t)}}|\nabla V|^2dx+\int_{B_{R_\ast(t)}}\left(\frac{n(n+1)}{r^2}+1\right)|V|^2dx\\
  &\leq C|\varepsilon|^2\|f_{n,m}^1\|^2_{L^{2}( B_{R_\ast(t)})}+C(n^2+1)|\widetilde{\rho}(t)|^2+C(n^2+1)|\varepsilon|^2|b^2_{n,m}(t)|^2.
\end{split}
\end{equation*}
Then integrating with $e^{2\delta_1t}dt$ and using (\ref{4.1}), (\ref{0.11}), (\ref{4.4}), we obtain
\begin{equation*}
\begin{split}
  \int_0^\infty e^{2\delta_1t}\|w-{w}_M\|^2_{H^{1}( B_{R_\ast(t)})}dt\leq& C|\varepsilon|^2\int_0^\infty e^{2\delta_1t}\|f^1\|^2_{L^{2}( B_{R_\ast(t)})}dt\\
  &+C\int_0^\infty e^{2\delta_1t}\|\rho-\widetilde{\rho}\|^2_{H^{7/2}(\partial B_{R_\ast(t)})}dt\\
  &+C|\varepsilon|^2\int_0^\infty e^{2\delta_1t}\|b^2\|^2_{H^{1}(\partial B_{R_\ast(t)})}dt\\
  \leq& C|\varepsilon|,
\end{split}
\end{equation*}
which completes (\ref{0.18}).

Using (\ref{0.18}) and the elliptic estimates, we derive (\ref{0.19}). It follows from (\ref{3.170}), (\ref{0.11}), (\ref{0.19}) and (\ref{4.3}) that (\ref{0.20}) holds. Therefore, the proof is complete.
 \hfill $\Box$

\section{Stability for $\mu<\mu_\ast$}
In this section, we shall show that the mapping $S$ defined in Remark \ref{rem1} admits a fixed point, thereby establishing the global existence of the system (\ref{1.2})-(\ref{1.6}) and the asymptotic stability of the radially symmetric $T$-periodic positive solution $(\sigma_\ast(r,t),p_\ast(r,t),R_\ast(t))$. Since the proof is similar to that of the corresponding proof in \cite[Page 633-637]{FH2}, we only give an outline of the proof and omit most of its details.

We first have the following H\"{o}lder estimates on $\rho$, $w$, $q$.

\begin{lem}\label{lemma9}
After performing a change of variables $x\rightarrow x+\varepsilon a^\ast(\varepsilon)$ on the initial data through (\ref{6.2}), then for all $t>0$, there exists a unique solution $(w,q,\rho)$ of the problem (\ref{3.150})-(\ref{3.190}) satisfying
\begin{equation}\label{0.30}
\begin{split}
  \|w\|_{C^{2+2\alpha/3,1+\alpha/3}(B_{R_\ast(t)}\times[0,\infty))}\leq C,
\end{split}
\end{equation}
\begin{equation}\label{0.31}
\begin{split}
  \|q\|_{C^{2+\alpha,\alpha/3}(B_{R_\ast(t)}\times[0,\infty))}\leq C,
\end{split}
\end{equation}
\begin{equation}\label{0.32}
\begin{split}
  \|\rho,D_x\rho\|_{C^{3+\alpha,1+\alpha/3}(\partial B_{R_\ast(t)}\times[0,\infty))}\leq C.
\end{split}
\end{equation}
\end{lem}
The proof of this lemma is similar to that of \cite[Lemma 7.1]{FH2}, so we omit here.

Introduce the space $X$ of functions $\Phi=(f^1,f^2,b^1,b^2,b^3)$, where the norm $\|\Phi\|$ defined by the maximum of the left-hand sides of (\ref{4.1})-(\ref{4.10}) with $\sqrt{|\varepsilon|}$ dropped, and set
\begin{equation*}
  X_1=\{\Phi\in X:\ \sqrt{|\varepsilon|}\|\Phi\|\leq1\},
\end{equation*}
then we define a new function $\widetilde{\Phi}\equiv S\Phi=(\widetilde{f}^1,\widetilde{f}^2,\widetilde{b}^1,\widetilde{b}^2,\widetilde{b}^3)$ as follows (cf. (\ref{0.17})):
\begin{equation*}
\begin{split}
  \widetilde{f}^1=A_\varepsilon w, \ \ \ \ \widetilde{f}^2=A_\varepsilon q,\ \ \ \ \ \widetilde{b}^1=B_\varepsilon^1,\ \ \ \ \ \widetilde{b}^2=B_\varepsilon^2,\ \ \ \ \ \widetilde{b}^3=B_\varepsilon^3,
\end{split}
\end{equation*}
where $w$, $q$, $\rho$ is the solution of (\ref{3.150})-(\ref{3.190}).
As in the proof of \cite[Lemma 7.2]{FH2} which is based on (\ref{0.30})-(\ref{0.32}) and (\ref{0.11})-(\ref{0.20}), we obtain that $S$ maps $X_1$ into itself and $S$ is a contraction mapping. Hence, the mapping $S$ admits a unique fixed point.

By the above argument, we can now state the main result of this paper.
\begin{thm}\label{th4}
Consider the problem (\ref{1.2})-(\ref{1.6}) with the initial data (\ref{0.68}) satisfying (\ref{0.57}) and let $\mu<\mu_\ast$. If $|\varepsilon|$ is sufficiently small, then there exists a unique global solution of this problem. And there exists a new center $\varepsilon a^\ast(\varepsilon)$, where $a^\ast(\varepsilon)$ is a bounded function of $\varepsilon$, such that for some $\overline{t}>0$, $\partial\Omega(t)$ behaves like
\begin{equation*}
  \partial B_{R_\ast(t)}(\varepsilon a^\ast(\varepsilon))=\big\{x:|x-\varepsilon a^\ast(\varepsilon)|=R_\ast(t)\big\}, \ \ \ \ \ t>\overline{t}.
\end{equation*}
\end{thm}

\begin{rem}
The proof of Theorem \ref{th4} shows that  after the translation of the origin $0\rightarrow\varepsilon a^\ast(\varepsilon)$ and the Hanzawa transformation, the global solution in the new variables $(r,\theta,\varphi)$ is of the form:
\begin{equation*}\label{3.1}
\begin{split}
    \sigma(r,\theta,\varphi,t)&=\sigma_\ast(r,t)+\varepsilon w(r,\theta,\varphi,t),\\
    p(r,\theta,\varphi,t)&=p_\ast(r,t)+\varepsilon q(r,\theta,\varphi,t),\\
    \partial\Omega(t):\ r&=R_\ast(t)+\varepsilon\rho(\theta,\varphi,t),
\end{split}
\end{equation*}
where $w$, $q$, $\rho$ satisfy (\ref{0.30})-(\ref{0.32}).
\end{rem}

\section{Acknowledgments}
This research is supported by China Postdoctoral Science Foundation (Grant No. 2020M683014). I would like to thank Professor S.B. Cui for his encouragement and helpful suggestions. I also thank Professor B. Hu for his advices in the preparation of this paper and constant encouragement.

\end{CJK}
\end{document}